\documentclass[a4, 12pt]{amsart}
\usepackage[mathscr]{eucal}
\usepackage{amssymb}
\usepackage{latexsym}
\usepackage{amsthm}
\theoremstyle{plain}
\newtheorem{theorem}{Theorem}[section]

\newtheorem{remark}{Remark}[section]
\newtheorem{lemma}{Lemma}[section]

\makeatletter
\@addtoreset{equation}{section}

\title[The second gap on  complete self-shrinkers]{The second gap on  complete self-shrinkers}
\author{Qing-Ming Cheng, Guoxin Wei and Wataru Yano}
\address{Qing-Ming Cheng \\ Department of Applied Mathematics, Faculty of Sciences,
Fukuoka  University, 814-0180, Fukuoka,  Japan, cheng@fukuoka-u.ac.jp}
\address{Guoxin Wei \\  School of Mathematical Sciences, South China Normal University,
510631, Guangzhou,  China, weiguoxin@tsinghua.org.cn}
\address{Wataru Yano \\ Department of Applied Mathematics, School  of Sciences,
Fukuoka  University, 814-0180, Fukuoka,  Japan,  kon.wata@gmail.com}
\begin{document}
\maketitle

\begin{abstract}
\noindent In this paper,  we study complete self-shrinkers in Euclidean space  and prove
 that an $n$-dimensional complete self-shrinker  in Euclidean space $\mathbb{R}^{n+1}$
 is  isometric to either  $\mathbb{R}^{n}$,  $S^{n}(\sqrt{n})$, or  $S^k (\sqrt{k})\times\mathbb{R}^{n-k}$, $1\leq k\leq n-1$,  if the squared norm $S$ of the second fundamental form, $f_3$ are constant and $S$ satisfies $S<1.83379$.  We should remark that  the condition of polynomial volume growth is not assumed.

\end{abstract}

\footnotetext{The first author was partially  supported by JSPS Grant-in-Aid for Scientific Research (B):  No.16H03937.
The second author was partly supported by grant No. 11771154 of NSFC, Guangdong Province Universities and Colleges Pearl River Scholar Funded Scheme (2018), Guangdong Natural Science Foundation Grant No.2019A1515011451.}

\section {Introduction}
\noindent

Let $X: M\rightarrow \mathbb{R}^{n+1}$ be a smooth $n$-dimensional immersed hypersurface in the $(n+1)$-dimensional
Euclidean space $\mathbb{R}^{n+1}$. One calls an immersed hypersurface $X: M\rightarrow \mathbb{R}^{n+1}$  a self-shrinker
 if it satisfies:
\begin{equation*}
\mathbf{H}=-X^{\bot},
\end{equation*}
where  $\mathbf{H}$ denotes the mean curvature vector of $M$, $X^{\bot}$ denotes the orthogonal
projection of $X$ onto the normal bundle of $M$.

It is well known that
Huisken \cite{H2}  \cite{H3} and  Colding and Minicozzi \cite{CM} have  proved that  if $M$ is  an $n$-dimensional  complete  embedded self-shrinker
in $\mathbb{R}^{n+1}$ with $H\geq 0$ and  with polynomial volume growth,
then $M$ is  isometric to either the  hyperplane  $\mathbb{R}^{n}$,
the  round sphere $S^{n}(\sqrt{n})$, or a cylinder $S^m (\sqrt{m})\times \mathbb{R}^{n-m}$, $1\leq m\leq n-1$.
For $n=1$, see Abresch and Langer \cite{AL}.

\begin{remark}
As one knows  that self-shrinkers play an important role in the study of the mean curvature flow because
they describe all possible blow up at a given singularity of a mean curvature flow.
\end{remark}

On the other hand,   Le and Sesum \cite{LS}  and  Cao and  Li \cite{CL}   have proved
that if $M$ is an $n$-dimensional complete self-shrinker  with polynomial volume growth
and $S\leq1$ in Euclidean space $\mathbb{R}^{n+1}$, then $M$ is isometric to either the hyperplane $\mathbb{R}^{n}$,  the  round sphere $S^{n}(\sqrt{n})$,  or a cylinder $S^m (\sqrt{m})\times \mathbb{R}^{n-m}$, $1\leq m\leq n-1$.
Ding and Xin \cite{DX1} have studied the second gap on the squared norm of the second fundamental form and they have
proved that if $M$ is  an $n$-dimensional complete self-shrinker with polynomial volume growth in Euclidean space $\mathbb{R}^{n+1}$,
there exists a positive number $\delta=0.022$ such that if $1\leq S\leq 1+0.022$, then $S=1$.
Furthermore, Cheng and Wei \cite{CW} have proved
\begin{theorem}
Let $X: M\rightarrow \mathbb{R}^{n+1}$
 be an $n$-dimensional complete self-shrinker  with polynomial volume growth in $\mathbb{R}^{n+1}$.
If the squared norm $S$ of the second fundamental form is constant and satisfies
$$S\leq 1+\frac{3}{7},$$
then $X: M\rightarrow \mathbb{R}^{n+1}$  is isometric to one of the following:
 \begin{enumerate}
\item $\mathbb{R}^{n}$,
\item
a cylinder $S^{k} (\sqrt k)\times \mathbb{R}^{n-k}$,
\item the round sphere $S^{n}(\sqrt{n})$.
\end{enumerate}
\end{theorem}

In \cite{CO}, Cheng and Ogata have obtained the following results (cf. Ding and Xin \cite{DX1}).
\begin{theorem}
Let $X: M\rightarrow \mathbb{R}^{3}$  be a $2$-dimensional complete self-shrinker in $\mathbb{R}^{3}$.
If the squared norm $S$ of the second fundamental form is constant,
 then $X: M\rightarrow \mathbb{R}^{3}$  is isometric to one of the following:
\begin{enumerate}
\item $\mathbb{R}^{2}$,
\item
a cylinder $S^1 (1)\times \mathbb{R}$,
\item the round sphere $S^{2}(\sqrt{2})$.
\end{enumerate}\end{theorem}

The following conjecture is known:

\vskip2mm
\noindent
{\bf Conjecture.}
Let $X: M\rightarrow \mathbb{R}^{n+1}$  be an $n$-dimensional complete self-shrinker in $\mathbb{R}^{n+1}$.
If the squared norm $S$ of the second fundamental form is constant,
 then $X: M\rightarrow \mathbb{R}^{n+1}$  is isometric to one of the following:
\begin{enumerate}
\item $\mathbb{R}^{n}$,
\item
a cylinder $S^{k} (\sqrt k)\times \mathbb{R}^{n-k}$,
\item the round sphere $S^{n}(\sqrt{n})$.
\end{enumerate}
\begin{remark} According to the result of Cheng and Ogata \cite{CO}, this conjecture has been   solved for $n=2$.
Recently, Cheng, Li and Wei \cite{clw} have solved this conjecture for $n=3$ under the condition that
$f_4$ is constant by making use of the generalized maximum principle due to Cheng and Peng \cite{CP}.
\end{remark}

For general $n$,  since this problem is too difficult, one can consider the special case:

\vskip2mm
\noindent
{\bf Problem.}
Let $X: M\rightarrow \mathbb{R}^{n+1}$  be an $n$-dimensional complete self-shrinker in $\mathbb{R}^{n+1}$.
If the squared norm $S$ of the second fundamental form is constant,
 then $X: M\rightarrow \mathbb{R}^{n+1}$  is isometric to one of the following:
\begin{enumerate}
\item $\mathbb{R}^{n}$,
\item
a cylinder $S^{k} (\sqrt k)\times \mathbb{R}^{n-k}$,
\item the round sphere $S^{n}(\sqrt{n})$,
\item $S\geq 2$.
\end{enumerate}

In this paper, we prove the following:
\begin{theorem}\label{theorem 1.3}
Let $X: M\rightarrow \mathbb{R}^{n+1}$  be an $n$-dimensional complete self-shrinker  in $\mathbb{R}^{n+1}$.
If the squared norm $S$ of the second fundamental form  and $f_3$ are  constants and $S$ satisfies
$$
S\leq 1.83379,
$$
then $M$ is isometric to one of the following:

\begin{enumerate}
\item $\mathbb{R}^{n}$,
\item
a cylinder $S^{k} (\sqrt k)\times \mathbb{R}^{n-k}$,
\item the round sphere $S^{n}(\sqrt{n})$,
\end{enumerate}
where $f_3=\sum_{j=1}^n\lambda_j^3$ and $\lambda_j$'s are principal curvatures of $X: M\rightarrow \mathbb{R}^{n+1}$.
\end{theorem}

\begin{remark} In our theorem \ref{theorem 1.3}, we do not assume that complete self-shrinkers $X: M\rightarrow \mathbb{R}^{n+1}$
have polynomial volume growth and it is known that there are many complete self-shrinkers without polynomial volume growth.
\end{remark}

\vskip 5mm

\section{Preliminaries}
\noindent
In this section, we give some notations and formulas.  Let $X: M\rightarrow \mathbb{R}^{n+1}$ be an $n$-dimensional self-shrinker in $\mathbb{R}^{n+1}$. Let $\{e_1,\cdots,e_{n},e_{n+1}\}$ be a local orthonormal basis along $M$ with dual
coframe $\{\omega_1,\cdots,\omega_{n},\omega_{n+1}\}$, such that $\{e_1,\cdots,e_{n}\}$ is a local
orthonormal basis of $M$ and $e_{n+1}$ is normal to $M$. Then we have
$$
\omega_{n+1}=0,\ \ \omega_{i n+1}=\sum_{j=1}^nh_{ij}\omega_j,\ \ h_{ij}=h_{ji},
$$
where $h_{ij}$ denotes the component of the second fundamental form of $M$.
$\mathbf{H}=\sum_{j=1}^nh_{jj}e_{n+1}$ is the mean curvature vector field, $H=|\mathbf{H}|=\sum_{j=1}^nh_{jj}$ is the mean curvature and $II=\sum_{i,j}h_{ij}\omega_i\otimes\omega_je_{n+1}$ is the second fundamental form of $M$.
The Gauss equations and Codazzi equations are given by
\begin{equation}\label{eq:12-6-4}
R_{ijkl}=h_{ik}h_{jl}-h_{il}h_{jk},
\end{equation}
\begin{equation}\label{eq:12-6-5}
h_{ijk}=h_{ikj},
\end{equation}
where $R_{ijkl}$ is the component of curvature tensor, the covariant derivative of $h_{ij}$ is defined by
$$
\sum_{k=1}^nh_{ijk}\omega_k=dh_{ij}+\sum_{k=1}^nh_{kj}\omega_{ki}+\sum_{k=1}^nh_{ik}\omega_{kj}.
$$
Let
$$F_i=\nabla_iF,\ F_{ij}=\nabla_j\nabla_iF,\ h_{ijk}=\nabla_kh_{ij}, \ {\rm and}\ h_{ijkl}=\nabla_l\nabla_kh_{ij},$$
 where
$\nabla_j $ is the covariant differentiation operator, we have
\begin{equation}\label{eq:12-6-6}
h_{ijkl}-h_{ijlk}=\sum_{m=1}^nh_{im}R_{mjkl}+\sum_{m=1}^nh_{mj}R_{mikl}.
\end{equation}
The following elliptic operator $\mathcal{L}$ is introduced by Colding and Minicozzi in \cite{CM}:
\begin{equation}\label{eq:12-6-1}
\mathcal{L}f=\Delta f-\langle X,\nabla f\rangle,
\end{equation}
where $\Delta$ and $\nabla$ denote the Laplacian and the gradient operator on the self-shrinker, respectively and $\langle\cdot,\cdot\rangle$ denotes the standard inner product of $\mathbb{R}^{n+1}$. By a direct calculation, we have
\begin{equation}\label{eq:12-6-2}
 \mathcal{L}h_{ij}=(1-S)h_{ij},\ \  \mathcal{L}H=H(1-S),\ \ \mathcal{L}X_i=-X_i,\ \ \mathcal{L}|X|^2=2(n-|X|^2),
\end{equation}
\begin{equation}\label{eq:12-6-3}
\frac{1}{2}\mathcal{L}S=\sum_{i,j,k}h_{ijk}^2+S(1-S).
\end{equation}
If $S$ is constant, then we obtain from  \eqref{eq:12-6-3}
\begin{equation}\label{eq:010}
\sum_{i,j,k}h_{ijk}^2=S(S-1),
\end{equation}
hence one has either
\begin{equation}
S=0,\ \ {\rm or}\ \ S=1,\ \ {\rm or}\ \ S>1.
\end{equation}
We can choose a local field of
orthonormal frames on $M^n$ such that, at the point that we consider,
$$
h_{ij}=\left \{\aligned  \lambda_i,\ \  \quad & \text {if} \quad i=j,\\
		    0,\ \  \quad  & \text {if} \quad i \neq j,\endaligned \right.
$$
then
$$
S=\sum_{i,j}h_{ij}^2=\sum_i\lambda_i^2,
$$
where $\lambda_i$ is called the principal curvature of $M$.
From \eqref{eq:12-6-4} and  \eqref{eq:12-6-6}, we get
\begin{equation}\label{eq:12-6-7}
h_{ijij}-h_{jiji}=(\lambda_i-\lambda_j)\lambda_i\lambda_j.
\end{equation}
By a direct calculation, we obtain
\begin{equation}
\sum_{i,j,k,l}h_{ijkl}^2=S(S-1)(S-2)+3(A-2B),
\end{equation}
where $A=\sum_{i,j,k}\lambda_i^2h_{ijk}^2,\ \ B=\sum_{i,j,k}\lambda_i\lambda_jh_{ijk}^2$.

We define two functions $f_3$ and $f_4$ as follows:
$$
f_3=\sum_{i,j,k} h_{ij}h_{jk}h_{ki}=\sum_{j=1}^n\lambda_j^3,\ \  \ \ f_4=\sum_{i,j,k,l} h_{ij}h_{jk}h_{kl}h_{li}=\sum_{j=1}^n\lambda_j^4,
$$
Then, the following formulas can be found in \cite{CW}:

\begin{lemma}
  Let $X: M\rightarrow \mathbb{R}^{n+1}$ be an $n$-dimensional  self-shrinker in
  $\mathbb{R}^{n+1}$. Then
  \begin{equation}\label{eq:12-6-8}
  \dfrac13\mathcal{L}f_3=(1-S)f_3+2C,
 \end{equation}
\begin{equation}\label{eq:12-6-9}
  \dfrac14\mathcal{L}f_4= (1-S)f_4+ (2A+B).
 \end{equation}
 \begin{equation}\label{eq:011} A-B \le \frac 13(\lambda_1-\lambda_2)^2tS^2,
\end{equation}
\begin{equation}
C^2\leq \dfrac13(A+2B)tS^2,
\end{equation}
where $C=\sum_{i,j,k} \lambda_ih_{ijk}^2$.
\end{lemma}

\vskip 5mm

 \section{Estimates for geometric invariants}
 \noindent
In this section, we will give some estimates which are needed to prove our theorem.
From now on, we denote
$$
S-1=tS,
$$
where $t$ is a positive constant and $t<\frac{1}{2}$ if we assume that $S$ is constant and $S>1$, then
$$(1-t)S=1, \ \ \ \ \sum_{i,j,k}h_{ijk}^2=tS^2.$$
Defining
$$
 u_{ijkl}:=\frac 14(h_{ijkl}+h_{jkli}+h_{klij}+h_{lijk} ),
 $$
we have
$$
  \sum_{i,j,k,l}h_{ijkl}^2 \geq \sum_{i,j,k,l}u_{ijkl}^2 +\frac 32  (Sf_4-f_3^2)
  $$
according to Gauss equations (2.1).\\
Let
\begin{equation}\label{eq:4-22-1}
Sf \equiv Sf_4-\bigl(f_3\bigl)^2=S \sum_i\lambda_i^4-\left(
\sum_i\lambda_i^3\right)^2=\dfrac12\sum(\lambda_i-\lambda_j)^2\lambda_i^2\lambda_j^2.
\end{equation}
Since $S$ and $f_3$ are constant, we have
$$
\sum_i\lambda_ih_{iik}=0, \ \  \sum_i\lambda_i^2h_{iik}=0, \ \text{\rm for any  } \ k
$$
and
$$
\sum_i\lambda_ih_{iikl}= -\sum_{i,j}h_{ijk}h_{ijl}, \ \  \sum_i\lambda_i^2h_{iikl}=
-2\sum_{i,j}\lambda_ih_{ijk}h_{ijl}, \ \text{\rm for any  } \ k, l.
$$
Hence, one has
$$
\sum_{i,j,}h_{iijj} \lambda_i\lambda_j=-C, \ \ \sum_{i,j}h_{iijj}\lambda_i^2\lambda_j=-2B.
$$
Defining
$$
a_{ij}=\sum_mh_{im}h_{mj}-yh_{ij},
$$
with $y=\dfrac{f_3}S$,
we have
$$
\sum_{i,j,k,l}u_{ijkl}h_{ij}h_{kl}=\sum_{i,j}\dfrac12(h_{iijj}+h_{jjii})\lambda_i\lambda_j=-C.
$$
$$
\sum_{i,j,k,l}u_{ijkl}a_{ij}h_{kl}=\sum_{i,j}\dfrac12(h_{iijj}+h_{jjii})(\lambda_{i}^2-y\lambda_i)\lambda_j=-B-\dfrac12A+yC.
$$
Hence, because of
$$
\sum_{i,j,k,l}\biggl\{u_{ijkl}+\alpha (a_{ij}h_{kl}+h_{ij}a_{kl})+\beta h_{ij}h_{kl}\biggl\}^2\geq 0,
$$
we obtain
\begin{equation*}
\begin{aligned}
\sum_{i,j,k,l}u_{ijkl}^2 &\geq -2\alpha \sum_{i,j,k,l}u_{ijkl}(a_{ij}h_{kl}+h_{ij}a_{kl})-\alpha^2\sum_{i,j,k,l}(a_{ij}h_{kl}+h_{ij}a_{kl})^2\\
&-2\beta \sum_{i,j,k,l}u_{ijkl}h_{ij}h_{kl}-\beta^2\sum_{i,j,k,l}(h_{ij}h_{kl})^2
-2\alpha\beta\sum_{i,j,k,l}(a_{ij}h_{kl}+h_{ij}a_{kl})h_{ij}h_{kl} \\
&=2\alpha(2B+A-2yC)-2\alpha^2Sf+2\beta C-\beta^2S^2 \\
&\geq 2\alpha(2B+A-2yC)-2\alpha^2Sf+\dfrac{C^2}{S^2}
\end{aligned}
\end{equation*}
by taking $\beta =\dfrac {C}{S^2}$.
Since $f_3$ is constant, we have from Lemma 2.1
$$
tSf_3=2C
$$
and
\begin{equation}\label{eq:4-28-1}
\begin{aligned}
Sf&=Sf_4-f_3^2=\dfrac12\sum_{i,j}(\lambda_i-\lambda_j)^2\lambda_i^2\lambda_j^2\\
&=\dfrac12\sum_{i,j}(h_{iijj}-h_{jjii})(\lambda_i-\lambda_j)\lambda_i\lambda_j\\
&=A-2B,
\end{aligned}
\end{equation}
that is,

\begin{equation}
Sf=Sf_4-f_3^2=A-2B.
\end{equation}
From
\begin{equation}
\sum_{i,j,k,l}h_{ijkl}^2=S(S-1)(S-2)+3(A-2B)
\end{equation}
and
$$
  \sum_{i,j,k,l}h_{ijkl}^2 =\sum_{i,j,k,l}u_{ijkl}^2 +\frac 32  (Sf_4-f_3^2),
  $$
  we have
\begin{equation}\label{eq:4-28-3}
S(S-1)(S-2)= \sum_{i,j,k,l}u_{ijkl}^2-\dfrac32(A-2B).
\end{equation}

Since $S$ is constant, by making use of  the generalized maximum  principle due to Cheng and Peng \cite{CP}, we have that there exist
a sequence $\{p_k\} $ in $M$ such that
\begin{equation}\label{eq:4-22-2}
tSf_4\geq 2A+B,
\end{equation}
it follows from \eqref{eq:4-28-1} that
$$
tf_3^2\geq  (2-t)A+(1+2t)B.
$$
Thus, for $z\geq 0$, from Lemma 2.1, one has
$$
C^2\leq  \dfrac13(A+2B)tS^2.
$$
\begin{equation}\label{eq:4-28-2}
\begin{aligned}
&\sum_{i,j,k,l}u_{ijkl}^2-\dfrac32(A-2B) \\
&\geq  2\alpha(2B+A-2yC)-2\alpha^2Sf+\dfrac{C^2}{S^2}-\dfrac32(A-2B)\\
&=2\alpha(2B+A)-(2\alpha^2+\dfrac32)(A-2B)+(-2\alpha t+(1+z)\dfrac{t^2}4)f^2_3-z\dfrac{C^2}{S^2}\\
&\geq 2\alpha(2B+A)-(2\alpha^2+\dfrac32)(A-2B)\\
&+(-2\alpha +(1+z)\dfrac{t }4)\bigl\{(2-t)A+(1+2t)B\bigl\}-z\dfrac{t}3(A+2B)\\
&=\bigl\{ 2\alpha -2\alpha^2 -\dfrac32+(-2\alpha +(1+z)\dfrac{t }4) (2-t) -z\dfrac{t}3\bigl\} A \\
&+\bigl\{ 4\alpha +4 \alpha^2+3  +(-2\alpha +(1+z)\dfrac{t }4)  (1+2t)  -z\dfrac{2t}3 \bigl\}B,\\
\end{aligned}
\end{equation}
where $-2 \alpha+(1+z)\frac{t}{4}\geq 0$.
By taking
$$
tz=\dfrac{8\alpha^2-8t\alpha +6+t(t+3)}{1-t},
$$
we have
\begin{equation}
\begin{aligned}
&4\alpha +4 \alpha^2+3  +(-2\alpha +(1+z)\dfrac{t }4)  (1+2t)  -z\dfrac{2t}3\\
&=-2\alpha +2\alpha^2 +\dfrac32-(-2\alpha +(1+z)\dfrac{t }4) (2-t) +z\dfrac{t}3\\
&=2(1-t)\alpha +2\alpha^2 +\dfrac32-\dfrac{t(2-t) }4 +tz\dfrac{3t-2}{12}\\
&=2(1-t)\alpha +2\alpha^2 +\dfrac32-\dfrac{t(2-t) }4 +\dfrac{8\alpha^2-8t\alpha +6+t(t+3)}{1-t}\dfrac{3t-2}{12}\\
&=\dfrac{4\alpha^2+4(3-4t)\alpha +3+2t(4t-3)}{6(1-t)}.
\end{aligned}
\end{equation}
We take $\alpha$ such that
\begin{equation}\label{eq:4-28-4}
4\alpha^2+4(3-4t)\alpha +3+2t(4t-3)=0
\end{equation}
if $t\leq \dfrac{9-\sqrt {33}}8$. Thus, we have from \eqref{eq:4-28-3}, \eqref{eq:4-28-2} and \eqref{eq:4-28-4} that $t\geq \dfrac12$. It is impossible.
Hence, we have
$$
t> \dfrac{9-\sqrt {33}}8.
$$
In this case,   taking $\alpha =-\dfrac{3-4t}2$, we obtain

$$
2\alpha^2+2(3-4t)\alpha +\dfrac32+t(4t-3)=-\dfrac{(3-4t)^2}2 +\dfrac32+t(4t-3)=-4t^2+9t-3,
$$
\begin{equation}\label{eq:4-22-11}
\begin{aligned}
tS^2(2t-1)S&\geq \sum_{i,j,k,l}u_{ijkl}^2-\dfrac32(A-2B) \\
&\geq \dfrac{4t^2-9t+3}{3(1-t)}( A -B).\\
\end{aligned}
\end{equation}

\noindent For any  $i , j$, we have
$$
 -\lambda_i\lambda_j\leq\dfrac14(\lambda_i-\lambda_j)^2.
$$
Hence, we get,  for $\lambda_i\lambda_j\leq 0$ and  $\lambda_i\lambda_k\leq 0$,
\begin{equation}\label{eq:12-23-4}
(|\lambda_i\lambda_j|+|\lambda_i\lambda_k|)^3\leq 4(|\lambda_i\lambda_j|^3+|\lambda_i\lambda_k|^3)
\leq (\lambda_i-\lambda_j)^2\lambda_i^2\lambda_j^2+(\lambda_i-\lambda_k)^2\lambda_i^2\lambda_k^2.
\end{equation}
For three different $i , j, k$, we know that at least one of
$\lambda_i\lambda_j$,  $\lambda_i\lambda_k$ and $\lambda_j\lambda_k$ is non-negative.
Without loss of generality, we assume
$\lambda_j\lambda_k\geq0$
and $\lambda_i\lambda_j\leq0$, $\lambda_i\lambda_k\leq0$,
then we get from \eqref{eq:12-23-4} and \eqref{eq:4-22-1}
\begin{equation}\label{eq:12-23-6}
-\lambda_i\lambda_j-\lambda_i\lambda_k-\lambda_j\lambda_k\leq |\lambda_i\lambda_j|
+|\lambda_i\lambda_k|
\leq (Sf)^{\frac13},
\end{equation}
and
\begin{equation}\label{eq:12-23-7}
-2\lambda_i\lambda_j\leq 2(|\lambda_i\lambda_j|^3)^{\frac13}
\leq 2(\dfrac14(\lambda_i-\lambda_j)^2\lambda_i^2\lambda_j^2)^{\frac13}
\leq 2\bigl(\frac{1}{4}Sf\bigl)^{\frac13}.
\end{equation}
Since
\begin{equation}
\begin{aligned}
3(A-B)&=
\sum_{i,j,k}(\lambda_i^2+\lambda_j^2+\lambda_k^2-\lambda_i\lambda_j-\lambda_i\lambda_k-\lambda_j\lambda_k)h_{ijk}^2\\
&=\sum_{i\neq j\neq k\neq i}(\lambda_i^2+\lambda_j^2+\lambda_k^2-\lambda_i\lambda_j-\lambda_i\lambda_k-\lambda_j\lambda_k)h_{ijk}^2\\
&\ \ \ +3\sum_{i\neq j}(\lambda_i^2+\lambda_j^2-2\lambda_i\lambda_j)h_{iij}^2,
\end{aligned}
\end{equation}
we conclude from \eqref{eq:12-23-6} and \eqref{eq:12-23-7},
\begin{equation}\label{eq:4-22-3}
\begin{aligned}
3(A-B)&\leq \biggl(S+2\biggl(\frac{1}{4}Sf\biggl)^{\frac13}\biggl)\sum_{i,j,k}h_{ijk}^2
=\biggl(S+2\biggl(\dfrac{Sf}4\biggl)^{\frac13}\biggl)tS^2.
\\
\end{aligned}
\end{equation}
From Lemma 2.1 and \eqref{eq:4-22-2}, one gets
\begin{equation}\label{eq:4-22-4}
Sf\leq \frac{A-B}{3(1-t)}.
\end{equation}
Since $S-1=tS$, we infer from \eqref{eq:4-22-3} and \eqref{eq:4-22-4}
\begin{equation}\label{eq:4-22-5}
\begin{aligned}
3(A-B)&\leq \biggl(S+2\biggl(\dfrac{A-B}{12(1-t)}\biggl)^{\frac13}\biggl)tS^2.
\end{aligned}
\end{equation}

\noindent
If we assume that $3(A-B)\leq a_k t S^3$, then we obtain from \eqref{eq:4-22-5} that
\begin{equation}
3(A-B)\leq a_{k+1}tS^3,
\end{equation}
where $$a_{k+1}= \biggl(1+2\biggl(\dfrac{a_k}{36}\frac{t}{1-t}\biggl)^{\frac13}\biggl).$$
Let $a_1=2$ and $t<0.454682$, we get from the above equations
\begin{equation}
a_7\leq 1.67738,
\end{equation}
then
\begin{equation}\label{eq:12-23-10}
3(A-B)\leq1.67738 t S^3.
\end{equation}
From \eqref{eq:4-22-11} and \eqref{eq:12-23-10}, we have
\begin{equation}\label{eq:12-23-11}
(2t-1)\geq \dfrac{4t^2-9t+3}{3(1-t)}\times \frac{1.67738 }{3}.
\end{equation}
Then, we get
\begin{equation*}
(t-0.454682)(t-1.24897)\leq0,
\end{equation*}
it follows that $t\geq 0.454682$. It is a contradiction. Hence, we get $t\geq0.454682$ and $S\geq 1.83379$.
We complete the proof of Theorem \ref{theorem 1.3}.

\end {document}